\documentclass[11pt]{article}

\usepackage{amssymb,amsmath,latexsym}
\usepackage[dvips]{graphics}
\usepackage{mathrsfs}

\newtheorem{thm}{Theorem}[section]

\newtheorem{cor}{Corollary}[section]

\newtheorem{lemma}{Lemma}[section]
\numberwithin{equation}{section}
\renewcommand{\baselinestretch}{1.4}
\def\ZZ{\mathbb{Z}}

\def\NN{\mathbb{N}}
\def\EE{\mathbb{E}}

\def\<{\langle}
\def\>{\rangle}
\def\pf{\noindent{\bf Proof.} }

\def\E{{\bf E}}

\def\P{{\bf P}}

\def\qed{{\hfill $\Box$\medskip}}

\allowdisplaybreaks[4]

\begin{document}
\title{\bf Some sufficient conditions for infinite collisions of simple random walks   on a wedge  comb}
\author{Xinxing Chen\footnote{Department of Mathematics, Shanghai Jiaotong University, Shanghai, China}~ and ~ Dayue Chen\footnote{School of Mathematical Sciences, Peking University, Beijing,
China }
 } \maketitle

\begin{abstract}
  In this paper, we give some sufficient conditions for the infinite collisions of
  independent simple random walks on a wedge comb with profile $\{f(n), n\in \ZZ\}$. One interesting
  result is that if $f(n)$ has a growth order as $n\log n$, then two independent simple random walks on the wedge comb will
  collide infinitely many times.
  Another is that if
  $\{f(n); n\in \ZZ\}$ are given by i.i.d. non-negative random  variables
  with finite mean, then for almost all  wedge comb with such profile,  three independent simple random walks on it will collide infinitely many times.

 \end{abstract} \noindent{\bf 2000 MR subject
classification:} 60K

\noindent {\bf Key words:} wedge comb, random walk, infinite
collision property, local time

\section{Introduction}
A  simple random walk on a graph  is defined as the Markov chain
that a particle  jumps from one vertex to a neighbor with equal
probability. Let $X=(X_n)$ and $X'=(X'_n)$ be two independent simple
random walks  starting from the same vertex.  We say that $X$ and
$X'$ collide infinitely often if $|\{n:X_n=X'_n\}|=\infty$.
 If   $X$ and $X'$ almost surely  collide infinitely often, then we say
that the graph has the infinite collision property. While if   $X$
and $X'$ almost surely collide finitely many times, then we say that
the graph has the finite collision property. Krishnapur and Peres
\cite{KP} first finds the example Comb$(\ZZ)$ on which  two simple
random random walks almost surely collide finitely many times, while
the expected number of collisions is infinite.

However, there is no simple monotonicity property for the finite
collision property, as exemplified by
$\ZZ\subset$Comb($\ZZ$)$\subset\ZZ^2$. Both $\ZZ$ and $\ZZ^2$ have
the infinite collision property but Comb($\ZZ$) has the finite
collision property. So, it is interesting to studying a subgraph of
Comb$(\ZZ)$.

{\bf Definition.}  Let $f$ be a function from $\ZZ$ to
$\mathbb{R}^+$. It induces a wedge comb Comb$(\ZZ,
f)=(\mathbb{V},\EE)$ which has the vertex set
$$
\mathbb{V}=\{(x,y):~x,y\in \ZZ, -f(x)\le y\le f(x)\}
$$
and  edge set $ \{[(x,n),(x,m)]: |m-n|=1\}\cup\{[(x,0),(y,0)]:
|x-y|=1\}$.

 Chen, Wei, and Zhang \cite{CWZ} shows
that Comb($\mathbb{Z}, f$)  has the infinite collision property when
$f(n)<n^\frac{1}{5}$. Recently, Barlow, Peres, Sousi \cite{BP} gives
a sufficient condition (in terms of Green functions) for infinite
collisions and   shows that Comb($\mathbb{Z}, f$) has the infinite
collision property when $f(n)\le n $; while it has the finite
collision property when $f(x)= n^\alpha$ for each $\alpha>1$.
Collisions on other graphes, such as random infinite cluster
 and random  tree,  can be seen in  \cite{BP}\cite{CCD}. In this paper, we focus only on the wedge combs
 with different profile $f(n)$,
and give a sufficient condition for a wedge comb which has the
infinite collision property, i.e.,
\begin{thm}\label{t:1.1} Let $
\breve{f}(n)=1\vee\max_{-n\le i\le n} f(i)$.  If
\begin{equation}\label{e:000001}
\sum_{n=1}^\infty\frac{1}{\breve{f}(n)}=\infty,
\end{equation}
then   two simple random walks on Comb($\ZZ,f$) will  collide
infinitely many times with probability one.  
\end{thm}

 As an directly application of Theorem \ref{t:1.1}, one has
\begin{cor}\label{c:1.2} If  $f(-n)+f(n)=O(n\log n)$, then two simple random walks on Comb($\ZZ,f$) will
collide  infinitely many times with probability one.
\end{cor}

Corollary \ref{c:1.2} improves the result of \cite{CWZ} and
\cite{BP}. On the other hand,  One can compare it with Theorem
\ref{t:4.1}, which says that for each $\beta>2$, if
$f(x)=|x|\log^\beta (|x|\vee 1)$ then Comb($\ZZ,f$) has the finite
collision property. But we guess that it should still have the
finite collision property for $1<\beta\le 2$.

A natural question to ask is what happen if there are more than two
independent simple random walks.  Suppose that $X$ is a simple
random walk and $X',X''$ are two independent copies of $X$. We say
that $X, X'$ and $X''$ collide together infinitely often if
$|\{n:X_n=X'_n=X''_n\}|=\infty$. it shows in \cite{BP} that three
independent simple random walks on $\ZZ$ will collide together
infinitely many times; while  on Comb($\ZZ,\alpha$) they will do
finitely many times for each $ \alpha>0$.
\begin{thm}\label{t:1.2}
Let $\{f(n); n\in \ZZ\}$ are independent and identically distributed
random variables with law $\mu$ supported in $[0,\infty)$. If $\mu$
has finite mean, then  for almost all $f$ three independent simple
random walks on Comb($\ZZ,f$) will collide  together infinitely many
times with probability one.
\end{thm}



\section{Proof of Theorem \ref{t:1.1}}
Let $X = \{X_n\}$ be a simple random walk on Comb$(\ZZ, f)$. Write
$$
X_n=(U_n, V_n)
$$
for all $n\ge 0$. So that  $U$ is a random process on $\ZZ$. For
each $k\ge 0$, we set $T_0=0$ and inductively
\begin{equation}\label{T_k/definition}
T_{k+1}:=\inf\{n>T_{k}: U_{n}\not= U_{n-1}\}.
\end{equation}
Then $X_n$ stays at the segment $\{(u,y):|y|\le f(u) \}$ when
$U_{T_{k}}=u$ during time $[T_{k}, T_{k+1}-1]$. Let
$$
W_k=U_{T_k}.
$$
It is easy to know such $\{W_k\}$ is a simple random walk on $\ZZ$
(by the strong Markov property).

 For any $a,b\in
\mathbb{R}$, we write $a\vee b=\max\{a,b\}$ and $a\wedge
b=\min\{a,b\}$. For a set $A$, let $|A|$ be the number of elements
of $A$. For each $x\in \mathbb{V}$, we write $x_1$ for the first
coordinate of $x$, $x_2$ for the second. For each $n$, let
$$
\mathbb{V}_n=\{(x_1, x_2)\in \mathbb{V}: |x_1|\le n\},
$$
and
$$
\theta_n=\inf\{m\ge 0: X_m\not\in \mathbb{V}_{n-1}\}.
$$
So, if $X_0\in \mathbb{V}_{n-1}$ then $\theta_n$ is the hitting time
of $\{(-n,0),(n,0)\}$ (i.e., the boundary of $\mathbb{V}_n$)  by
$X$.

Let $X'$ be another simple random walk on Comb$(\ZZ, f)$,
independent of $X$. Define
$U_n', V_n', W'_k, T_k',\theta_{n}'$ as well. For each pair $u, v\in
\mathbb{V}$,  We write $\P^u$ for the probability measure of a
simple random walk $X$ starting from $u$ and  write $\P^{u,v}$ for
the joint probability measure of the two independent simple random
walks $X$ and $X'$ starting from $u$ and $v$, respectively. We
 also write $\E^u$ and $\E^{u,v}$ for the corresponding expectations.
For each $m\ge 0$, set  $\sigma_0=0$, and inductively
$$
\sigma_{m+1}:=\inf\{n> \sigma_m: ~U_n=U_n', ~  U_n\not=U_{n-1}
{\rm~or~} U_n'\not=U_{n-1}'\}.
$$

\begin{lemma}\label{l:2.1} For any $\varepsilon>0$, there exist integers $d$ and $N_0$, such that
for all $N\ge N_0$ and all $u,v\in \mathbb{V}_N$ with
$u_1+u_2+v_1+v_2$ even,
$$
\P^{u,v}\left(\sigma_{N}\ge \theta_{dN}\wedge\theta_{dN}'\right)<
\varepsilon.
$$
\end{lemma}
\pf Fix $\varepsilon>0$. Suppose that $N\in \NN$ is large enough and
that $u,v\in \mathbb{V}_N$. For each $n\ge 0$, set
$$
Z_{2n}=U_n-U_n'~~~~{\rm and~~~~}Z_{2n+1}=U_{n+1}-U_n'.
$$
Set $\tau_0=0$, and inductively
$$
\tau_m=\inf\{n>\tau_{m-1}: Z_n \not=Z_{n-1}\}.
$$
By the strong Markov property, $\{Z_{\tau_m}, m\ge 0\}$ is a simple
random walk on
  $\ZZ$  and
\begin{equation}\label{e:002}
Z_{\tau_0}=u_1-v_1\in[-2N, 2N].
\end{equation}
 If $Z_{\tau_m}=0$ then
  \begin{align*}
  &U_{k}=U_{k}',~U_k\not=U_{k-1}',~\tau_m=2k~~{\rm for~~some~} k\in
  \ZZ^+; {\rm~~or~}\\
  & U_{k+1}=U_k',~U_{k}\not=U_k',~\tau_m=2k+1~~{\rm for~~some~} k\in
  \ZZ^+.
  \end{align*}
Notice that
  $U_n+V_n+U_n'+V_n'$ is always even under the assumption that $U_0+V_0+U_0'+V_0'=u_1+u_2+v_1+v_2$ is even.
This fact, together with $U_{k+1}=U_k'\not=U_k$, implies that
$U_{k+1} =U_{k+1}'$. For each $M\ge 0$, let
$$
\xi({0, M}):=|\{m: Z_{\tau_m}=0, 0\le m\le  M\}|,
$$
the local time of 0 by $\{Z_{\tau_m}, m\ge 0\}$. As a result of the
previous argument,
\begin{equation}\label{e:1.1}
\{\xi(0,M)\ge N\}\subseteq\{\sigma_N \le \tau_M\}\subseteq\{\sigma_N
\le T_M\wedge T_M'\}.
\end{equation} By  (9.11) on Page 39 of  \cite{PR} and (\ref{e:002}), there exists $x\in \NN$ such that
\begin{equation}\label{e:001}
\P^{u,v}(\xi(0,x N^2)\le  N)< \frac{\varepsilon}{2}.
\end{equation}
On the other hand, by Theorem 2.13  on Page 21 of \cite{PR}, we can
find such $d\in \NN$ which satisfies the following inequality.
\begin{align*}
& \P^{u,v}(T_{x N^2}\ge\theta_{dN} {\rm ~or~} T_{x
N^2}'\ge\theta_{dN}') \\
= & \P^{u,v}\left(\max_{0\le k\le x N^2} |W_k|\ge dN{\rm
~or~}\max_{0\le k\le x N^2} |W_k'|\ge dN\right) \le
\frac{\varepsilon}{2}.
\end{align*}
Together with (\ref{e:1.1}) and (\ref{e:002}), we can get the
desired result. \qed

The above lemma shows that, with a small exception, the number of
collisions is bounded by departures times $\theta$ and $\theta '$ of
$U$ and $U'$ linearly.

Secondly, we estimate the probability of that there is at least one
collision of $X$ and $X'$ once $U$ and $U'$ collide. For every $m\ge
0$, define
\begin{align*}
\Psi_m=\big\{   X_n=X'_n,~ &|V_n|+|V_n'|\ge
|V_{\sigma_m}|+|V_{\sigma_m}'|~\\
 &{\rm~for~some~}\sigma_m\le
n<\inf\{h> \sigma_m: V_h=0 {\rm ~or~}V_h'=0 \}\big\}.
\end{align*}
Notice that if $X_{\sigma_m}=X_{\sigma_m}'$ then $\Psi_m$ occurs.
Moreover for every $(u,0), (u,v)\in \mathbb{V}$ with ${v}$ being
even, if $\Psi_m\cap \{X_{\sigma_m}=(u,0), X_{\sigma_m}'=(u,v)\}$
occurs then $X$ enters the segment $L=\{(u,x): 1\le x\le f(u)\}$ and
collides with $X_n'$ at a height not less than ${v}/{2}$ after time
$\sigma_m$ but before one of them leaving $L$. Here, by height we
mean the second component of a vertex $x\in \mathbb{V}$.
We need these events in order to have good bounds as follows.

\begin{lemma}\label{l:2.2} There exist positive numbers $c_1$ and $c_2$ such that for all $ (u,v)\in \mathbb{V}$ with $v$ being even,
$$
\frac{c_1}{|v|\vee 1}\le \P^{(u,0), (u,v)}( \Psi_0)\le
\frac{c_2}{|v|\vee1}.
$$
\end{lemma}

\pf First let us examine the case that $f(u)\ge  2v$  and $v$ is
even. For each $x\in \ZZ$, define
$$
\tau_x=\inf\{n>0: X_n=(u, x)\}.
$$
Define $\tau_x'$ similarly. If $\tau_{2v}\le \tau_0$ and $X'$ stay
in $[{v}/{2}, {3v}/{2}]$ before time $\tau_{2v}$, then $X$ and $X'$
must collide before $\tau_{2v}\le \tau_0\wedge \tau_0'$ at a vertex
whose height is greater than or equal to ${v}/{2}$. Therefore
\begin{align*}
\P^{(u,0),(u,v)}(\Psi_0)\ge&\P^{(u,0),(u,v)}\left(X_1=(u,1),\
\tau_{2v}\le v^2\wedge \tau_0,\  \tau_{v/2}'\vee \tau_{3v/2}'\ge
v^2\right)~~~~~~~~~~~~~~~~~~~~~~~~~~\\
=&\P^{(u,0)}(X_1=(u,1))\P^{(u,1)}\left(\tau_{2v}<  v^2, \tau_{2v}\le
\tau_0\right)\P^{(u,v)}\left(\tau_{v/2}\vee \tau_{3v/2}\ge
v^2\right)
\end{align*}
\begin{equation}\label{e:006}
~~~\ge \frac{1}{4}\P^{(u,1)}\left(\tau_{2v}< v^2, \tau_{2v}\le
\tau_0\right)\P^{(u,v)}\left(\tau_{v/2}\vee \tau_{3v/2}\ge
v^2\right).
\end{equation}
Once $X$ enters into the  segment $L$, by observation,  the behavior
of $V$, the second component of $X$, before $V$ hitting $0$, has the
same law as a simple random walk on $\ZZ\cap [0, f(u)]$. Let
$\{\eta_i\}$ be i.i.d. random variables with
$$
\P(\eta_1=1)=\P(\eta_1=-1)=1/2,
$$
and
$$
\widetilde{\tau}_{2v}=\inf\left\{k>0: 1+\sum_{i=1}^k\eta_i\ge
2v\right\}, ~~\widetilde{\tau}_{0}=\inf\left\{k>0:
1+\sum_{i=1}^k\eta_i\le 0\right\}.
$$
Then
\begin{equation}
\P^{(u,1)}\left(\tau_{2v}< v^2, \tau_{2v}\le
\tau_0\right)=\P\left(\widetilde{\tau}_{2v}< v^2,
\widetilde{\tau}_{2v}\le \widetilde{\tau}_0\right).
\end{equation}
Obviously, $\{\widetilde{\tau}_{2v}\le  v^2\}$ and
$\{\widetilde{\tau}_{2v}\le \widetilde{\tau}_0\}$ are both
increasing event.  By the FKG inequality,
\begin{equation}\label{e:005}
\P\left(\widetilde{\tau}_{2v}\le  v^2, \widetilde{\tau}_{2v}\le
\widetilde{\tau}_0\right)\ge \P(\widetilde{\tau}_{2v}\le
v^2)\P(\widetilde{\tau}_{2v}\le \widetilde{\tau}_0).
\end{equation}
By Lemma 3.1 of \cite{PR}, 
\begin{equation}\label{e:003}
\P(\widetilde{\tau}_{2v}\le \widetilde{\tau}_0)= \frac{1}{2v}.
\end{equation}
By Theorem 2.13 of \cite{PR} again, there exists $c_1>0$
independently of $u$ and $v$, such that
\begin{equation}\label{e:004}
\P(\widetilde{\tau}_{2v}\le v^2)\ge c_1 {\rm~~~~~ and~~~~~}
\P^{(u,v)}\left(\tau_{v/2}\vee \tau_{3v/2}\ge v^2\right)\ge c_1.
\end{equation}
Taking (\ref{e:006})-(\ref{e:004}) together,  we obtain the first
inequality of the lemma.
\begin{equation}
\P^{(u,0),(u,v)}(\Psi_0)\ge \frac{c_1^2}{8v}. \end{equation}
\\

Now we turn to proving the second inequality. Define
$$
H=\sum_{n=0}^{\infty}1_{\{V_n=V_n',~ n< \tau_0\wedge \tau_0'\}},
$$
the number of collisions of $X$ and $X'$ before  one of them leaving
$L$. Then
\begin{align*}
& \E^{(u,0),(u,v)}(H)=
\sum_{n=0}^{\infty}\sum_{x=1}^{f(u)}\P^{(u,0),(u,v)}(V_n=V_n'=x,~
n\le \tau_0\wedge
\tau_0')\\
= &\sum_{n=0}^{\infty}\sum_{x=1}^{f(u)}\P^{(u,0)}(V_n=x,~ n\le
\tau_0)\P^{(u,v)}(V_n=x,~ n\le
\tau_0)\\
\le & 2\sum_{n=0}^{\infty}\sum_{x=1}^{f(u)}\P^{(u,0)}(X_n=(u,x),~
n\le \tau_0)\P^{(u,x)}(V_n=v,~ n\le
\tau_0)\\
\le & 2\sum_{n=0}^{\infty}\P^{(u,0)}(V_{2n}=v,~ 2n\le
\tau_0)\\
= & 2\E^{(u,0)}({\rm ~number~of~ visits~ to ~} (u,v)
{\rm~by~}X~{\rm~before~returning~to~} (u,0)~).
\end{align*}
In the previous arguments, the first inequality follows by the
knowledge of reversible Markov chain that for all $x\in [1,
f(u)]\cap \ZZ$
$$
\P^{(u,v)}(V_n=x,~ n\le \tau_0)\le 2 \P^{(u,x)}(V_n=v,~ n\le
\tau_0).
$$
By Theorem 9.7 of \cite{PR},
\begin{equation}\label{e:0010}
\E^{(u,0),(u,v)}(H)\le 2.
\end{equation}
The second inequality will follow once we show that there exists
$c_2>0$ independent of $u,v$ such that
\begin{equation}\label{e:009}\E^{u,v}(H| \Psi_0)\ge c_2v.
\end{equation}
Since we condition on the event $\Psi_0$, there is a collision at
position $x=(u, w)$ for some $w$ with $ w\ge {v}/{2}$. Conditioned
on this event, the total number of collisions that happen in the set
$\{(u,h): h\ge {v}/{3}\}$, will be greater than the number of
collisions that take place before the first time that one of the
random walks exits this interval. The lower bound could be obtained
by the following consideration.

Consider two independent simple random walks in an interval,
starting at $v/2$. Before hitting either $v/3$ or $2v/3$, the
average number of collisions is the number of average number of
returning to the starting point before exiting the interval. The
average number of returning to the starting point is exactly the
Green function of a simple random walk, starting at $v/2$, before
exiting the interval $(v/3, 2v/3)$, which is of order $v$.

By (\ref{e:0010}) and (\ref{e:009}), we have
$$
\P^{(u,0),(u,v)}(\Psi_0)\le \frac{2}{c_2v}.
$$
This completes the proof of the case that $f(u)\ge 2v$ and $v$ is
even. The proof can be modified to treat other cases and is omitted
here. \qed
\\
\\

By Lemma \ref{l:2.1}, we can find  $d\in \NN$ and $N_0\in \NN$ such
that for all $N\ge N_0$ and all $u,v\in \mathbb{V}_N$ with
$u_1+u_2+v_1+v_2$ even,
\begin{equation}\label{e:0020}
\P^{u,v}\left(\sigma_{N}\ge \theta_{dN}\wedge\theta_{dN}'\right)<
\frac{1}{2}.
\end{equation}
Fix $d$ through this section. To be concise, we set
$$
\breve{f}(n)=1\vee\max_{-n\le i\le n}f(i).
$$
As a result,  $\breve{f}(n)$ is a strictly positive and  increasing
function on $\ZZ^+$.

\begin{lemma}\label{l:000007} There exist $N_0\in \NN$ and $c>0$, such that for all $N\ge
N_0$ and all $u,v\in \mathbb{V}_{N}$ with $u_1+u_2+v_1+v_2$ even,
$$
\P^{u,v}\big(X_n=X_n' {\rm ~for~ some~}n\in [0,\theta_{dN}\wedge
\theta_{dN}')\big)\ge\frac{cN}{\breve{f}(dN)+N}.
$$
\end{lemma}
\pf Let
$$
H=\sum_{m=1}^N
\big(|V_{\sigma_m}|\vee|V_{\sigma_m}'|\vee1\big)1_{\Psi_m}1_{\{\sigma_m<
\theta_{dN}\wedge\theta_{dN}'\}}.
$$
As a result of that,  if $H>0$ then $X$ and $X'$ collide before they
break out of $\mathbb{V}_{dN}$. We shall use the second moment
method to estimate the probability of the occurrence of  $\{H>0\}$.
\begin{align*}
&~~~\E^{u,v}(H)\\
&=\sum_{m=1}^N\E^{u,v}\left((|V_{\sigma_m}|\vee|V_{\sigma_m}'|\vee1)1_{\Psi_m}1_{\{\sigma_m<
\theta_{dN}\wedge\theta_{dN}'\}}\right)\\
&=\sum_{m=1}^N\E^{u,v}\left(\E^{u,v}\left((|V_{\sigma_m}|\vee|V_{\sigma_m}'|\vee1)1_{\Psi_m}1_{\{\sigma_m<
\theta_{dN}\wedge\theta_{dN}'\}}~\big|~X_i, X_i', 0\le i\le \sigma_m\right)\right)\\
&=\sum_{m=1}^N\E^{u,v}\left((|V_{\sigma_m}|\vee|V_{\sigma_m}'|\vee1)\P^{u,v}\left(\Psi_m~\big|~X_i,
X_i', 0\le i\le \sigma_m\right); \sigma_m<
\theta_{dN}\wedge\theta_{dN}'\right)\\
&=\sum_{m=1}^N\E^{u,v}\left((|V_{\sigma_m}|\vee|V_{\sigma_m}'|\vee1)\P^{X_{\sigma_m},X_{\sigma_m}'}\left(\Psi_0\right);
\sigma_m<
\theta_{dN}\wedge\theta_{dN}'\right)\\
&\ge\sum_{m=1}^N\E^{u,v}\left(c_1; \sigma_m<
\theta_{dN}\wedge\theta_{dN}'\right)\\
&\ge c_1N\P^{u,v}( \sigma_N< \theta_{dN}\wedge\theta_{dN}')\\
&\ge \frac{c_1}{2}N,
\end{align*}
where the last inequality is by (\ref{e:0020}); the last three
inequality is by Lemma \ref{l:2.2} and the last  equation is by the
strong Markov property. Using  Lemma \ref{l:2.2} and the strong
Markov property again, we have
\begin{align*}
&\E^{u,v}(H^2)\\
=&\E^{u,v}\left(\sum_{m=1}^N
(|V_{\sigma_m}|\vee|V_{\sigma_m}'|\vee1)1_{\Psi_m}1_{\{\sigma_m<
\theta_{dN}\wedge\theta_{dN}'\}}\sum_{n=1}^N
(|V_{\sigma_n}|\vee|V_{\sigma_n}'|\vee1)1_{\Psi_n}1_{\{\sigma_n<
\theta_{dN}\wedge\theta_{dN}'\}}    \right)\\
=&\E^{u,v}\left(\sum_{m=1}^N
(|V_{\sigma_m}|\vee|V_{\sigma_m}'|\vee1)^21_{\Psi_m}1_{\{\sigma_m<
\theta_{dN}\wedge\theta_{dN}'\}}\right)\\
&~~~+2\E^{u,v}\left(\sum_{m=1}^N
(|V_{\sigma_m}|\vee|V_{\sigma_m}'|\vee1)1_{\Psi_m}1_{\{\sigma_m<
\theta_{dN}\wedge\theta_{dN}'\}}\sum_{n>m}^N
(|V_{\sigma_n}|\vee|V_{\sigma_n}'|\vee1)1_{\Psi_n}1_{\{\sigma_n<
\theta_{dN}\wedge\theta_{dN}'\}}    \right)\\
\le& 2\breve{f}(dN)\E^{u,v}\left(\sum_{m=1}^N
(|V_{\sigma_m}|\vee|V_{\sigma_m}'|\vee1)1_{\Psi_m}1_{\{\sigma_m<
\theta_{dN}\wedge\theta_{dN}'\}}\right)\\
&~~~~+2\sum_{m=1}^N\sum_{n>m}^N\E^{u,v}\left(
(|V_{\sigma_m}|\vee|V_{\sigma_m}'|\vee1)1_{\Psi_m}1_{\{\sigma_m<
\theta_{dN}\wedge\theta_{dN}'\}}
(|V_{\sigma_n}|\vee|V_{\sigma_n}'|\vee1)1_{\Psi_n}1_{\{\sigma_n<
\theta_{dN}\wedge\theta_{dN}'\}}    \right)\\
=&2\breve{f}(dN)\E^{u,v}(H)+2\sum_{m=1}^N\sum_{n>m}^N\E^{u,v}\big((|V_{\sigma_m}|\vee|V_{\sigma_m}'|\vee1)1_{\Psi_m}1_{\{\sigma_m<
\theta_{dN}\wedge\theta_{dN}'\}}\\
&~~~~~~~~~~~~~~~~~~~~~~~~~~~~~~~~~~~~~~~~~~~~~~~\times~(|V_{\sigma_n}|\vee|V_{\sigma_n}'|\vee1)
\P^{X_{\sigma_n},X_{\sigma_n}'}(\Psi_0)1_{\{\sigma_n<
\theta_{dN}\wedge\theta_{dN}'\}}\big)\\
\le
&2\breve{f}(dN)\E^{u,v}(H)+2c_2\sum_{m=1}^N\sum_{n>m}^N\E^{u,v}\left((|V_{\sigma_m}|\vee|V_{\sigma_m}'|\vee1)1_{\Psi_m}1_{\{\sigma_m<
\theta_{dN}\wedge\theta_{dN}'\}}\right)\\
\le& (2\breve{f}(dN)+2c_2N)\E^{u,v}(H).
\end{align*}
So that by the H$\ddot{o}$lder inequality,
$$
\P^{u,v}(H>0)\ge \frac{[\E^{u,v}(H)]^2}{\E^{u,v}[H^2]}\ge
\frac{\E^{u,v}(H)}{2\breve{f}(dN)+2c_2N}\ge
\frac{c_1N}{4\breve{f}(dN)+4c_2N}.
$$
Hence we have the result.\qed
\\
\\
\\
{\it Proof of Theorem \ref{t:1.1}.} We need only to prove the case
that two independent simple random walks on Comb($\ZZ,f$) starting
from the same vertex $(0,0)$. So, simply write
$\P=\P^{(0,0),(0,0)}$. For each $m\ge 1$, define
$$
\Upsilon_m=\{X_n=X_n' {\rm~~~~for~~some~~} n\in
[\theta_{d^m}\wedge\theta_{d^{m}}',\theta_{d^{m+1}}\wedge\theta_{d^{m+1}}')\}.
$$
Notice that for all $0\le n\le \theta_{d^m}\wedge\theta_{d^{m}}'$,
$$
X_n, X_n'\in \mathbb{V}_{d^m}.
$$
So, by the strong Markov property and Lemma \ref{l:000007}, there
exists $c>0$ such that for all $m$ large enough
\begin{align*}
&\P(\Upsilon_m|~1_{\Upsilon_i}, 1\le i <m, ~X_n,X_n', ~n\le
\theta_{d^m}\wedge\theta_{d^{m}}')\\
=&\P^{X_{t},X_{t}'}\big(X_n=X_n' {\rm~~~~for~~some~~} n\in
[0,\theta_{d^{m+1}}\wedge\theta_{d^{m+1}}')\big)\\
\ge& \frac{c d^{m+1}}{\breve{f}(d^{m+1})+d^{m+1} },
\end{align*}
where $ t=\theta_{d^m}\wedge\theta_{d^{m}}'.$ We shall show later
that (\ref{e:000001}) implies
\begin{equation}\label{e:000002}
\sum_{m=1}^\infty\frac{d^{m}}{\breve{f}(d^{m})+d^{m} }=\infty.
\end{equation}
If (\ref{e:000002}) holds, then by the second Borel Cantelli Lemma
(extend version, Page 237 of \cite{DR}),
$$
\P(~\Upsilon_m {\rm~~~infinitely~~~often}~)=1.
$$
Furthermore,
$$
\P(X_n=X_n'{\rm~~~infinitely~~~often}~)\ge  \P(~\Upsilon_m
{\rm~~~infinitely~~~often}~)=1.
$$

Now we prove that (\ref{e:000002}) holds. If $\breve{f}(d^m)\le d^m$
for infinitely many $m$, then (\ref{e:000002}) holds obviously.
Otherwise, there exists $m_0$ such that for all $m\ge m_0$,
$$
\breve{f}(d^m)> d^m.
$$
Hence to prove  (\ref{e:000002}), we need only to prove
\begin{equation}
\sum_{m=1}^\infty \frac{d^m}{\breve{f}(d^m)}=\infty.
\end{equation}
Set $S_m=d+d^2+\cdots+d^m$. Then
$$
   d^m\le S_m\le d^{m+1}.
$$
As a result,
$$
\frac{d^{m}}{\breve{f}(d^{m})}\ge \frac{1}{d}
\sum_{l=S_{m}+1}^{S_{m+1}}\frac{1}{\breve{f}(l)}.
$$
So
$$
\sum_{m=1}^\infty \frac{d^{m}}{\breve{f}(d^{m})}\ge
\frac{1}{d}\sum_{m=1}^\infty
\sum_{l=S_{m}+1}^{S_{m+1}}\frac{1}{\breve{f}(l)}\ge
\frac{1}{d}\sum_{l=d+1}^\infty\frac{1}{\breve{f}(l)}=\infty.
$$
 Such we completed the proof of
Theorem \ref{t:1.1}. \qed

\section{Proof of Theorem \ref{t:1.2}}

Theorem \ref{t:1.2} will follow once we show the following result.

\begin{thm} \label{t:3:1}If
$$\sum_{i=-n}^n f(i)=O(n),$$
then three independent simple random walks on Comb($\ZZ,f$) will
collide together infinitely many times with probability one.
\end{thm}

Let $X''$ be another  independent simple random walk on Comb($\ZZ,
f$). 
For each  $u, v,w\in \mathbb{V}$, we write $\P^{u,v,w}$ for the
joint probability measure of the three independent simple random
walks $X,X'$ and $X''$ starting from $u,v$ and $w$, respectively. We
write $\E^{u,v,w}$ for the corresponding expectation.

By the condition of Theorem \ref{t:3:1}, there exists $c>2$, such
that for all $n\in \ZZ^+$, \begin{equation}\label{e:3.1}
\sum_{i=-n}^n f(i)\le \left(\frac{c}{2}-2\right)n.
\end{equation}
We fix $f$ and $c$ which satisfy (\ref{e:3.1}) through this section.
Since
$$
[-n, n]\times\{0\}\cap \ZZ^2\subseteq
\mathbb{V}_n\subseteq\bigcup_{i=-n}^n [-f(i),f(i)]\times\{i\},
$$
one has
\begin{equation}\label{e:3.2}
2n\le |\mathbb{V}_n|\le cn.
\end{equation}
 For each $x\in \mathbb{V}$, let
$$
\tau_x=\inf\{m\ge 0:~X_m=x\},
$$
the hitting time of $x$ by $X$. Similarly, we have he following
estimates.

\begin{lemma}\label{l:3.1} For any $c_1>0$, there exist  $d\in \NN$, $c_2\in \NN$ and $n_0\in \NN$ such that,
for all $n\ge n_0$ and all  $u,v\in \mathbb{V}_n$,
$$
\P^u(\theta_{dn}\le n^2)\le c_1,
$$
and
$$\P^u\big( \tau_v> c_2 n^2 \big)\le c_1,$$
\end{lemma}
\pf  Suppose $n\in \NN$ is large enough and is fixed. Fix $u,v\in
\mathbb{V}_n$ and $c_1 \in R$. Assume further that $0< c_1 < 1/5$.

The first statement is obvious.  By Theorem 2.13 of \cite{PR},
 there exists $d\in\NN$, $d> 2$ and is independent of $n$ and $u$,
such that
\begin{equation}\label{e:3.25}
\P^u(\theta_{dn}\le n^2)\le \P^u\left( \max_{0\le k\le n^2}|W_k|\ge
dn\right)\le c_1.
\end{equation}

We now prove the second statement. By the strong Markov property,
$$
\P^u(\tau_v\le (c_2+c_4) n^2)\ge \P^u(\tau_{\widetilde{u}}\le c_2
n^2)\P^{\widetilde{u}}(\tau_v\le  T_{c_3n^2},~T_{c_3n^2}\le  c_4 n^2
).
$$
 To prove the second result, we need
only to prove that there exist $c_2,c_3\in \NN$  and $c_4>0$ which
are independent of $n$ and $u$  and satisfy
(\ref{e:3.8})-(\ref{e:3.10}) as follows. Here
$\widetilde{u}=(u_1,0)$ and $T_k$ is defined in
(\ref{T_k/definition}).
\begin{equation}\label{e:3.8}
\P^u(\tau_{\widetilde{u}}> c_2 n^2)\le c_1;
\end{equation}
\begin{equation}\label{e:3.9}
\P^{\widetilde{u}}(\tau_v> T_{c_3 n^2})\le c_1;
\end{equation}
\begin{equation}\label{e:3.10}
\P^{\widetilde{u}}(T_{c_3 n^2}> c_4 n^2)\le 3c_1.
\end{equation}
Once we verify that (\ref{e:3.8})-(\ref{e:3.10}) hold, then
$$
\P^u(\tau_v\le (c_2+c_4) n^2)\ge 1-5c_1.
$$

 Now we prove (\ref{e:3.8})-(\ref{e:3.10}) one by one. First, since $u\in \mathbb{V}_n$,
$$
f(u_1)\le \frac{cn}{2}.
$$
The process $X$  starting from $u$ stays at the segment $\{(u_1,y):
|y|\le f(u_1)\}$ before reaching $\widetilde{u}$. So its behavior
before reaching $\widetilde{u}$ is much like a simple random walk on
$\ZZ$. As a result of that, we can find such  $c_2$ as the
requirement of (\ref{e:3.8}).

Next,
 for each $x\in \ZZ$ and $n\in
\ZZ^+$, let
$$
\xi(x,n)=|\{k: W_k=x, 0\le k\le n\}|
$$
the  local time of $x$ before time $n$ by $W$. By  Theorem 9.4 and
(9.11)  of \cite{PR}, for any $c_1^*, c_2^*>0$ there exist
$c_3\in\NN$ independent of $u, v$ and $n$ but depending on $c_1^*$
and $c_2^*$, such that
\begin{equation}\label{e:3.001}
\P^{\widetilde{u}}(\xi(v_1,c_3n^2)\ge c_2^*n)\ge 1- c_1^*.
\end{equation}
By a similar argument as Lemma \ref{l:2.2}, we can get
\begin{equation}\label{e:3.002}
\P^{(v_1,0)}( X_n=v {\rm~ for~ some~} n\in [0, T_1))\ge
\frac{1}{4f(v_1)+1}\ge \frac{1}{4cn}.
\end{equation}
If we let $c_1^*$ small enough and $c_2^*$ large enough, then by
(\ref{e:3.001}),(\ref{e:3.002}) and  a proof as Theorem \ref{t:1.1},
we can find $c_3$ for (\ref{e:3.9}).

Finally, we estimate $T_n$.  Let
$$
\xi(n)=\max_{x\in \ZZ} \xi(x,n).
$$
According to Theorem 9.14 and  (10.6) of \cite{PR}, there exists
$c_1^*>0$ independent of $u$  and $n$ such that,
\begin{equation}\label{e:3.3}
\P^{\widetilde{u}}(\xi(c_3 n^2)> c_1^*n)\le c_1.
\end{equation}
For each $x\in \ZZ$ and $ k\in \NN$, let
$$D^x_k=\inf\big\{T_l: T_l> \inf\{n: \xi(x,n)=k\}   \big\} -\inf\{n: \xi(x,n)=k\}.$$
Hence $D_k^x$ is the time spent on the line segment $\{(x, y):
|y|\le f(x)\}$ at the $k$-th visit $x$ by $W$. Obviously, $\{D_k^x,
k\ge 1\}$ are independent and identically distributed. Moreover,
\begin{equation}\label{e:3.12}
\E^{\widetilde{u}}(D_k^x)=\E^{(x_1,0)}(D_1^x)\le 2f(x)+1.
\end{equation}
Condition on $V_0=0$, for all $n$ there has
\begin{equation}\label{e:3.13}
 T_{n}\le\sum_{x\in
\ZZ}\sum_{k=1}^{\xi({x, n})} D_k^x.
\end{equation}
 By (\ref{e:3.1}), (\ref{e:3.12}) and  (\ref{e:3.12}),
\begin{align*}
&\E^{\widetilde{u}}(T_{c_3n^2}; \theta_{c_3dn} \ge c_3 n^2, \xi(c_3n^2)\le c_1^*n)\\
\le&\E^{\widetilde{u}}\left(\sum_{x\in \ZZ}\sum_{k=1}^{\xi({x, c_3n^2})} D_k^x ;~\theta_{c_3dn} \ge c_3n^2, \xi(c_3n^2)\le c_1^*n\right)\\
\le&\E^{\widetilde{u}}\left(\sum_{x=-c_3dn}^{c_3dn}\sum_{k=1}^{c_1^* n} D_k^x ;~\theta_{c_3dn} \ge c_3n^2, \xi(c_3n^2)\le c_1^*n\right)\\
\le&\sum_{x=-c_3dn}^{c_3dn}\sum_{k=1}^{c_1^* n}\E^{\widetilde{u}}(D_k^x)\le c_1^*n\sum_{x=-c_3dn}^{c_3dn} (2f(x)+1)\\
\le& 3c_1^*c_3cdn^2
\end{align*}
By the Markov inequality
$$
\P^{\widetilde{u}}(T_{c_3n^2}> 3c_1^{-1}c_1^*c_3cd n^2;
\theta_{c_3dn} \ge c_3 n^2, \xi(c_3n^2)\le c_1^*n)\le c_1.
$$
This, together with (\ref{e:3.25}) and (\ref{e:3.3}), verifies
(\ref{e:3.10}). We have completed the proof.\qed
\\

By Lemma \ref{l:3.1}, we can find $c_2,d\in \NN$, such that  there
exists $n_0\in \NN$, for all  $n\ge n_0$ and all $u,x\in
\mathbb{V}_n$,
\begin{equation}\label{e:3.52}
\P^{u}(\tau_x> c_2 n^2)\le \frac{1}{4}
{\rm~~~and~~~}\P^{u}(\theta_{dn}\le c_2 n^2)\le \frac{1}{4}.
\end{equation}

\begin{lemma}\label{l:3.2} There exist $c_5>0$ and $N_0\in \NN$, such that for any integer  $N> N_0$ and all  $u,v,w\in \mathbb{V}_N$ with
$(-1)^{u_1+u_2}=(-1)^{v_1+v_2}=(-1)^{w_1+w_2}$,
$$
\P^{u,v,w}\left(X_n=X_n'=X_n'' {\rm~ for~ some~ }n\in \big[0,~
\theta_{dN}\wedge\theta_{dN}'\wedge\theta_{dN}''\big)\right)\ge
\frac{c_5}{\log N}.
$$
\end{lemma}
\pf Fix $N\in \NN$. For conciseness,  we write
$$
\Theta=\theta_{dN}\wedge\theta_{dN}'\wedge\theta_{dN}''=\inf\big\{m\ge
0: \{X_m,X_m',X_m''\}\not\subseteq \mathbb{V}_{dN-1}\big\}.
$$
Then  $\Theta$ is the first time that one of  $X^*$ breaks out of
$\mathbb{V}_{dN-1}$.
 Let
$$
H=\sum_{n=0}^{2c_2 N^2}1_{\{X_n=X_n'=X_n''\in \mathbb{V}_{N},~
\Theta> n\}}.
$$
We need only to prove \begin{equation}\label{e:3.41}
\E^{u,v,w}(H)\ge c_1^* {\rm~~~~~and ~~~~} \E^{u,v,w}(H^2)\le
c_2^*\E^{u,v,w}(H)\log N.
\end{equation}
Then  by H$\ddot{o}$lder inequality,
$$
\P^{u,v,w}(X_n=X_n'=X_n'' {\rm~ for~ some~ }n\leq \Theta)\ge
\P^{u,v,w}(H>0)\ge
\frac{\big[\E^{u,v,w}(H)\big]^2}{\E^{u,v,w}\big(H^2\big)}\ge
\frac{c_1^*}{c_2^*\log N}.
$$

Now we prove (\ref{e:3.41}). Let
$$
q_n(u,x)=\P^u(X_{n}=x, \theta_{dN}>n).
$$
Then $q_{2n}(x,x)$ is decreasing in $n$ for every $x\in
\mathbb{V}_N$ (Refer to \cite{BP}). By the strong Markov property,
for each $u,x\in \mathbb{V}_N$ and $c_2N^2\le n\le 2c_2 N^2$ with
$u_1+u_2+x_1+x_2+n$ even
\begin{align*}
\P^u(X_n=x, \theta_{dN}> n)
=&\P^u(X_n=x, \tau_x\le n < \theta_{dN})\\
=&\E^u( \P^x(X_{n-k}=x, \theta_{dN}>n)|_{\tau_x=k}; \tau_x\le n,
\theta_{dN}>\tau_x)\\
=&\E^u( q_{n-\tau_x}(x,x); \tau_x\le n,
\theta_{dN}>\tau_x)\\
\ge& q_{2c_2 N^2}(x,x)\P^u(\tau_x\le n, \theta_{dN}>\tau_x)\\
\ge&q_{2c_2 N^2}(x,x)\big(\P^u(\tau_x\le c_2 N^2)+\P^u(\theta_{dN}>
c_2 N^2)-1\big).
\end{align*}
By (\ref{e:3.52}),
$$
\P^u(X_n=x, \theta_{dN}> n)\ge \frac{1}{2}q_{2c_2 N^2}(x,x).
$$
As a result,
\begin{align*}
\E^{u,v,w}(H) & = \sum_{x\in
\mathbb{V}_{N}}\sum_{n=0}^{2c_2N^2}\P^{u,v,w}(X_n=X_n'=X_n''=x,\Theta> n)\\
&=\sum_{x\in \mathbb{V}_{N}}\sum_{n=0}^{2c_2N^2}\P^u(X_n=x,
\theta_{dN}> n)\P^v(X_n=x, \theta_{dN}> n)\P^w(X_n=x,
\theta_{dN}> n)\\
&\ge \frac{1}{8}\sum_{x\in
\mathbb{V}_{N}}~~\sum_{c_2N^2\le n\le 2c_2 N^2\  ~u_1+u_2+x_1+x_2+n {\rm ~even}~}\left[q_{2c_2 N^2}(x,x)\right]^3\\
 &\ge \frac{c_2N^2}{16}\sum_{x\in
\mathbb{V}_{N}}\left[q_{2c_2 N^2}(x,x)\right]^3\\
&\ge\frac{c_2N^2}{16|\mathbb{V}_N|^2}\left[\sum_{x\in
\mathbb{V}_{N}}q_{2c_2 N^2}(x,x)\right]^3 \ge
\frac{c_2}{16c^2}\left[\sum_{x\in \mathbb{V}_{N}}q_{2c_2
N^2}(x,x)\right]^3.
\end{align*}
Where the last inequality is by the H$\ddot{o}$lder inequality.
Using the H$\ddot{o}$lder inequality and  (\ref{e:3.52}) again, we
have
\begin{align*}
\sum_{x\in \mathbb{V}_{N}}q_{2c_2N^2}(x,x)\ge& \sum_{x\in
\mathbb{V}_{N}}\sum_{y\in
\mathbb{V}_{dN}}q_{c_2N^2}(x,y)q_{c_2N^2}(y,x)\\
\ge&\frac{1}{4}\sum_{x\in
\mathbb{V}_{N}}\sum_{y\in \mathbb{V}_{dN}}[q_{c_2N^2}(x,y)]^2\\
\ge&\frac{1}{4|\mathbb{V}_{dN}|}\sum_{x\in
\mathbb{V}_{N}}\left(\sum_{y\in
\mathbb{V}_{dN}}q_{c_2N^2}(x,y)\right)^2\\
\ge&\frac{|\mathbb{V}_{N}|}{4|\mathbb{V}_{dN}|}\min_{x\in\mathbb{V}_{N}}\left[\sum_y \P^x(X_{c_2N^2}=y, \theta_{dN}>c_2N^2)\right]^2\\
\ge
&\frac{1}{2cd}\min_{x\in\mathbb{V}_{N}}\left[\P^x(\theta_{dN}>c_2N^2)\right]^2
\ge\frac{1}{4cd}.
\end{align*}
Take together,
\begin{align*}
\E^{u,v,w}(H)\ge \frac{c_2}{16c^2}\cdot
\left(\frac{1}{4cd}\right)^3\ge \frac{c_2}{2000c^5 d^3}.
\end{align*}
So we have gotten  the first part of (\ref{e:3.41}).

Now we turn to the second moment, i.e., the second part of
(\ref{e:3.41}). Since that Comb($\ZZ, f$) is a graph with uniformly
bounded degree, there exists $c_1^*>0$, such that for all $x,y\in
\mathbb{V}$
\begin{align*}
\P^x(X_k=y)\le \frac{c_1^*}{\sqrt{k}}.
\end{align*}
Hence,
\begin{align*}
\P^{x,x,x}\left( X_{k}=X_{k}'=X_{k}''\right) =&\sum_{y\in
\mathbb{V}}\E^{x,x,x}\big(X_{k}=X_{k}'=X_{k}''=y\big)
\\=\sum_{y\in \mathbb{V}}\big[\P^{x}(X_k=y)\big]^3\le&
\frac{(c_1^*)^2}{k}\sum_{y\in \mathbb{V}}\P^{x}(X_k=y)
=\frac{(c_1^*)^2}{k}.
\end{align*}
 By the inequality above and  the strong Markov property,
\begin{align*}
&\E^{u,v,w}\left(H^2\right)\\
=&\E^{u,v,w}\left(\sum_{n=0}^{2c_2N^2}\sum_{x\in
\mathbb{V}_{N}}1_{\{X_n=X_n'=X_n''=x,\Theta>
n\}}\sum_{k=0}^{2c_2N^2} \sum_{y\in
\mathbb{V}_{N}}1_{\{X_k=X_k'=X_k''=y,\Theta>
k\}}\right)\\
\le&\E^{u,v,w}(H)+2\E^{u,v,w}\left(\sum_{n=0}^{2c_2N^2}\sum_{x\in
\mathbb{V}_{N}}1_{\{X_n=X_n'=X_n''=x,\Theta>
n\}}\sum_{k>n}^{2c_2N^2} 1_{\{X_k=X_k'=X_k''\}}\right)\\
=&\E^{u,v,w}(H)+2\sum_{n=0}^{2c_2N^2}\sum_{x\in
\mathbb{V}_{N}}\sum_{k>n}^{2c_2N^2}\P^{u,v,w}\left(X_n=X_n'=X_n''=x,\Theta>
n, X_k=X_k'=X_k''\right)\\
=&\E^{u,v,w}(H)+2\sum_{n=0}^{2c_2N^2}\sum_{x\in
\mathbb{V}_{N}}\sum_{k>n}^{2c_2N^2}\P^{u,v,w}\left(X_n=X_n'=X_n''=x,\Theta>
n\right)\P^{x,x,x}\left( X_{k-n}=X_{k-n}'=X_{k-n}''\right)\\
\le &\E^{u,v,w}(H)+2\sum_{n=0}^{2c_2N^2}\sum_{x\in
\mathbb{V}_{N}}\P^{u,v,w}\left(X_n=X_n'=X_n''=x,\Theta>
n\right)\sum_{k=1}^{2c_2N^2}\P^{x,x,x}\left( X_{k}=X_{k}'=X_{k}''\right)\\
\le & \E^{u,v,w}(H)\left(1+2\max_{x\in
\mathbb{V}_N}\sum_{k=1}^{2c_2N^2}\P^{x,x,x}\left(
X_{k}=X_{k}'=X_{k}''\right)\right)\\
\le&\E^{u,v,w}(H)\left(1+2\max_{x\in
\mathbb{V}_N}\sum_{k=1}^{2c_2N^2}\frac{(c_1^*)^2}{k}\right)\\
\le&\left(1+2(c_1^*)^2\log (2c_2)+4(c_1^*)^2\log N
\right)\E^{u,v,w}(H).
\end{align*}
Therefore, we have verified the second part of (\ref{e:3.41}) and
finished the proof of the lemma. \qed
\\
\\
{\it Proof of Theorem \ref{t:3:1}.}  For each $m\ge 1$, define event
$$
\Upsilon_m=\{X_n=X_n'=X_n'' {\rm~~~~for~~some~~} n\in [\theta_{d^m
}\wedge \theta_{d^m }' \wedge \theta_{d^m }'',\theta_{d^{m+1}
}\wedge \theta_{d^{m+1} }' \wedge \theta_{d^{m+1} }'')\}.
$$
Then by the strong Markov property and Lemma \ref{l:3.2}, for all
$m$ large enough
\begin{align*}
&\P(\Upsilon_m|1_{\Upsilon_i}, 1\le i<m, ~X_n,X_n',X_n'', ~n\le
\Theta_{d^m})\\
=&\P^{X_{t},X_{t}',X_t''}\big(X_n=X_n'=X_n'' {\rm~for~some~} n\in
[\theta_{d^m }\wedge \theta_{d^m }' \wedge \theta_{d^m
}'',\theta_{d^{m+1}
}\wedge \theta_{d^{m+1} }' \wedge \theta_{d^{m+1} }'')\big)\\
\ge& \frac{c_5}{\log d^m }=\frac{c_5}{m\log d},
\end{align*}
where $ t=\theta_{d^m}\wedge\theta_{d^{m}}'\wedge \theta_{d^m }''.$
By the second Borel Cantelli Lemma (extend version, Page 237 of
\cite{DR}),
$$
\P(~\Upsilon_m {\rm~~~infinitely~~~often}~)=1.
$$
Furthermore,
$$
\P(X_n=X_n'=X_n''{\rm~~~infinitely~~~often}~)\ge  \P(~\Upsilon_m
{\rm~~~infinitely~~~often}~)=1.
$$
Thus we completed the proof of Theorem \ref{t:3:1}. \qed

\section{A Related Result }
\begin{thm}\label{t:4.1}
Let $f(x)=|x|\log^\beta (|x|\vee1)$ for all $x\in \ZZ$. If
$\beta>2$, then the total number of collisions by
 two independent simple random walks on Comb($\ZZ,f$) is almost surely
finite.
\end{thm}
  The proof of the theorem   is almost the same as the case
$f(x)=x^\alpha$ for $\alpha>1$ in \cite{BP}.  So we just outline the
changes needed to run the proof.
\begin{lemma}\label{l:4.1} Let $f(x)=|x|\log^\beta (|x|\vee1)$ for all $x\in \ZZ$. Let $x=(k,h)\in \mathbb{V}$. Then the transition density $q$
satisfies:
$$
q_t(0,x)\le \frac{c}{n^2\log^\beta n}          {\rm
~~if~~}t=n^3\log^\beta n {\rm~~and~~}n\ge k,
$$
$$
q_t(0,x)\le \frac{c}{k^2\log^\beta k}          {\rm
~~if~~}t=n^3\log^\beta n {\rm~~and~~}n< k.
$$
\end{lemma}
\pf Prove similarly as  Lemma 5.1 of \cite{BP}.\qed
\\
 Set $Q_{k,h}$, where $h\le k\log^\beta k$, as follows:
$$
Q_{k,h}=\{(k,y): 0\le y\le h\}.
$$
We set $Z_{k,h}=Z(Q_{k,h})$ to be the number of collisions of the
two random walks in $Q_{k,h}$. We also define
$\widetilde{Z}=Z_{k,2h/3}-Z_{k,h/3}$, i.e. the number of collisions
that happen in the set $\{(k,y): \frac{h}{3}\le y\le
\frac{2h}{3}\}$.
\begin{lemma}\label{l:4.2}
$\E(Z_{k,h})\le {ch}/({k\log^\beta k})$;
and $\E(Z_{k,h}|\widetilde{Z}_{k,h}>0)\ge ch$. 
\end{lemma}

\pf By Lemma \ref{l:4.1},
\begin{align*}
\E(Z_{k,h})=&\sum_t\sum_{x\in Q_{k,h}}
q_t(0,x)^2\le\sum_{t<k^3\log^\beta k} \frac{ch}{k^4\log^{2\beta}
k}+\sum_{t\ge k^3\log^\beta k} \frac{ch}{g(t)^2}\\
\le&\frac{ch}{k\log^\beta k}+ch\int_{k^3\log^\beta
k}^\infty\frac{1}{g(t)^2}dt\\
=&\frac{ch}{k\log^\beta k}+ch\int_{k}^\infty\frac{1}{h^4\log^{2\beta} h}d(h^3\log^\beta h)\\
=&\frac{ch}{k\log^\beta k}+ch\int_{k}^\infty\frac{3}{h^2\log^{\beta}
h}dh
+ch\int_{k}^\infty\frac{\beta}{h^2\log^{\beta+1} h}dh\\
\le&\frac{c'h}{k\log^\beta k},
\end{align*}
where $n=g(t)$ is the inverse function of $t=n^3\log^\beta n$. The
second inequality is proved similarly as Lemma 5.2 of \cite{BP}.
\qed
\\
\\
{\it Proof of Theorem \ref{t:4.1}.} By Lemma \ref{l:4.2}
$$
\P(\widetilde{Z}_{k,h}>0)\le \frac{1}{k\log^\beta k}.
$$
Now summing over all $k$ and over all $h$ ranging over powers of 2
and satisfying $h\le k\log^\beta k$, we get that
$$
\sum_k \sum_{h {\rm~power~of~} 2} \P(\widetilde{Z}_{k,h}>0)\le
\sum_{k} \frac{\log_2 (k\log^\beta k) }{k\log^\beta k}<\infty,
{\rm~since~}\beta>2.
$$
Hence the total number of collisions is finite almost surely.\qed



\begin{thebibliography}{20}\small
\bibitem{BP} Barlow, M.T., Peres, Y., Sousi, P., Collisions of
Random Walks, {\em Preprint}, (2010).
\bibitem{CCD} Chen, X., Chen, D., Two random walks on the open cluster of $\ZZ^2$ meet infinitely
often. {\em Science China Mathematics}, {\bf 53}, 1971-1978 (2010).
\bibitem{CWZ} Chen. D., Wei, B. and Zhang, F., A note  on the finite
collision property of random walks. {\em  Statistics and Probability
Letters}, {\bf 78}, 1742-1747, (2008).
\bibitem{DR} Durrett, R.  {\em
Probability: Theory and Examples}, 3rd ed. Brooks/Cole, Belmont,
2005.
\bibitem{KP} Krishnapur, M. and Peres,Y., Recurrent graphs where two
independent random walks collide finitely often. {\em Elect. Comm.
in Probab.} {\bf 9}, 72-81, (2004).
\bibitem{Polya}Polya, G.,  {\em  George Polya: Collected Papers}, Volume IV,
582-585, The MIT Press, Cambridge, Massachusetts.
\bibitem{PR} R$\acute{e}$v$\acute{e}$sz, P.,  {\em Random walk in
random and non-random environments.} World Scientific Publishing Co.
Pte. Ltd.,   2005.
\renewcommand{\baselinestretch}{.1}
\end{thebibliography}
\end{document}